\documentstyle[12pt,psfig]{article}
\title{Some planar isospectral domains}
\author{Peter Buser \and John Conway \and Peter Doyle
\and Klaus-Dieter Semmler}
\date{Version 1.0.1\\15 September 1994}
\newcommand{\fig}[3]{
\begin{figure}
\psfig{figure=figures/#1.ps,width=370pt}
\caption{#3}
\label{#2}
\end{figure}
}
\newcommand{\by}{\ /\ }
\newcommand{\By}{\\&/& }
\newcommand{\ditto}{\prime \prime}
\newcommand{\cross}{\times}
\newcommand{\0}{\bf 0}
\newcommand{\1}{\bf 1}
\newcommand{\2}{\bf 2}
\newcommand{\3}{\bf 3}
\newcommand{\4}{\bf 4}
\newcommand{\5}{\bf 5}
\newcommand{\6}{\bf 6}
\begin{document}
\maketitle
\begin{abstract}
We give a number of examples of isospectral pairs of plane domains,
and a particularly simple method of proving isospectrality.
One of our examples
is a pair of domains that are not only isospectral but {\em homophonic}:
Each domain has a distinguished point such that
corresponding normalized Dirichlet eigenfunctions
take equal values at the distinguished points.
This shows that one really can't hear the shape of a drum.
\end{abstract}

\section{Introduction}
In 1965, Mark Kac
\cite{kac:drum}
asked, `Can one hear the shape of a drum?', so
popularizing the question of whether there can exist two non-congruent
isospectral domains in the plane.  In the ensuing 25 years many
examples of isospectral manifolds were found, whose dimensions,
topology, and curvature properties gradually approached those of
the plane.  Recently, Gordon, Webb, and Wolpert
\cite{gordonWebbWolpert:drum}
finally reduced
the examples into the plane.  In this note, we give a number of
examples, and a particularly simple method of proof.
One of our examples
(see Figure \ref{homophonic})
is a pair of domains that are not only isospectral but {\em homophonic}:
Each domain has a distinguished point such that
corresponding normalized Dirichlet eigenfunctions
take equal values at the distinguished points.
We interpret this to mean that
if the corresponding `drums'
are struck at these special points, then they
`sound the same' in the very strong sense that every frequency will
be excited to the same intensity for each.
This shows that one really can't hear the shape of a drum.

\fig{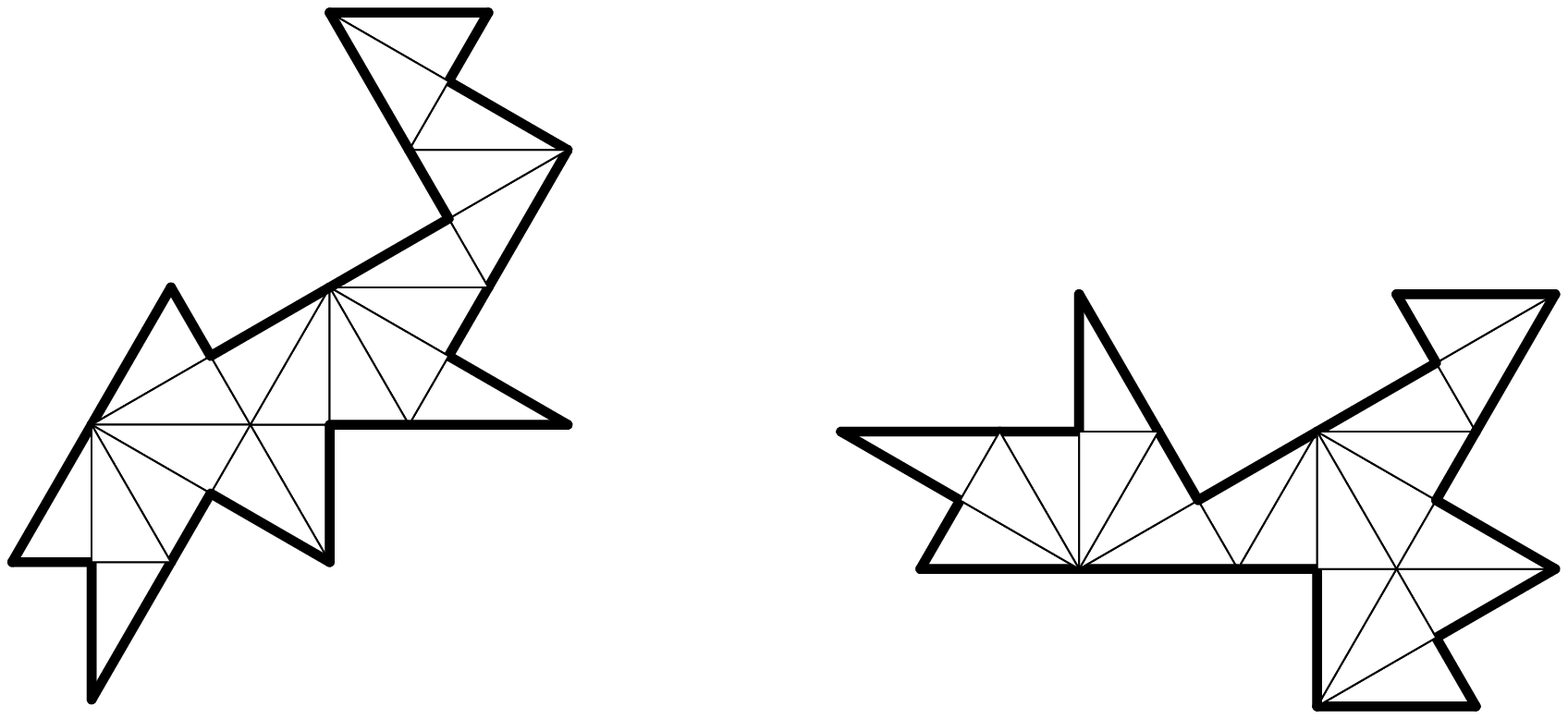}{homophonic}{Homophonic domains.  These drums sound the
same when struck at the interior points
where six triangles meet.}

\section{Transplantation}

The following transplantation proof was first applied to
Riemann surfaces by Buser \cite{buser:riemann}.
For our domains this proof turns out
to be particularly easy.

\fig{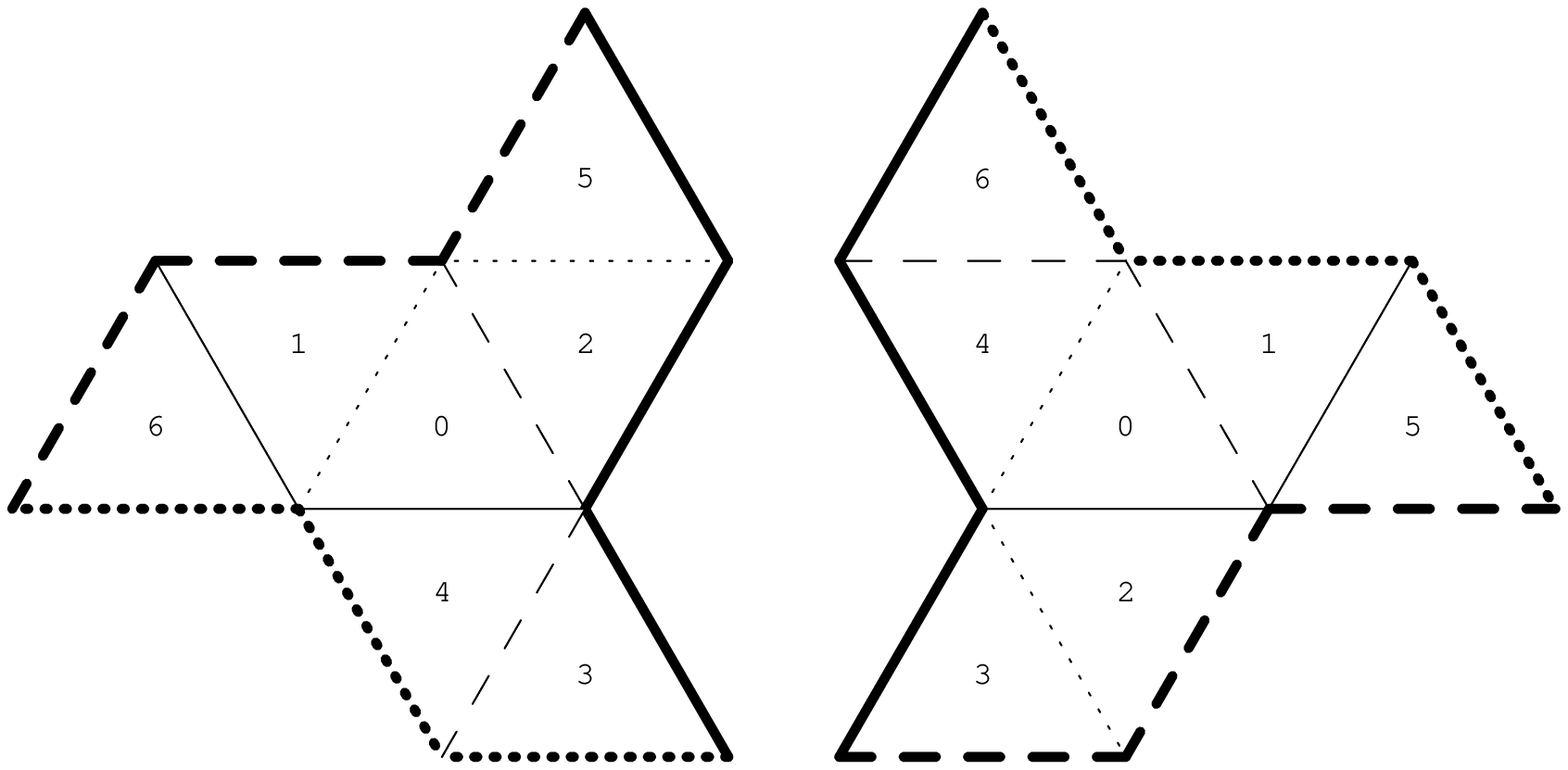}{propeller}{Propeller example.}

Consider the two propeller-shaped regions shown in Figure \ref{propeller}.
Each region consists of seven equilateral triangles (labelled in some
unspecified way).
Our first pair of examples is obtained from these by replacing the
equilateral triangles by acute-angled scalene triangles,
all congruent to each other. The propellers are triangulated by
these triangles in such a way that any two triangles
that meet along a line are mirror
images in that line, as in Figure \ref{warped}.
In both propellers the central triangle has a distinguishing property:
its sides connect the three inward corners of the propeller. The position of
the propellers in Figure \ref{warped} is such
that the unique isometry from the central triangle on the left-hand
side to the central triangle on the right-hand side is a translation.
This translation does not map the propellers onto one another and so
they are not isometric.

\fig{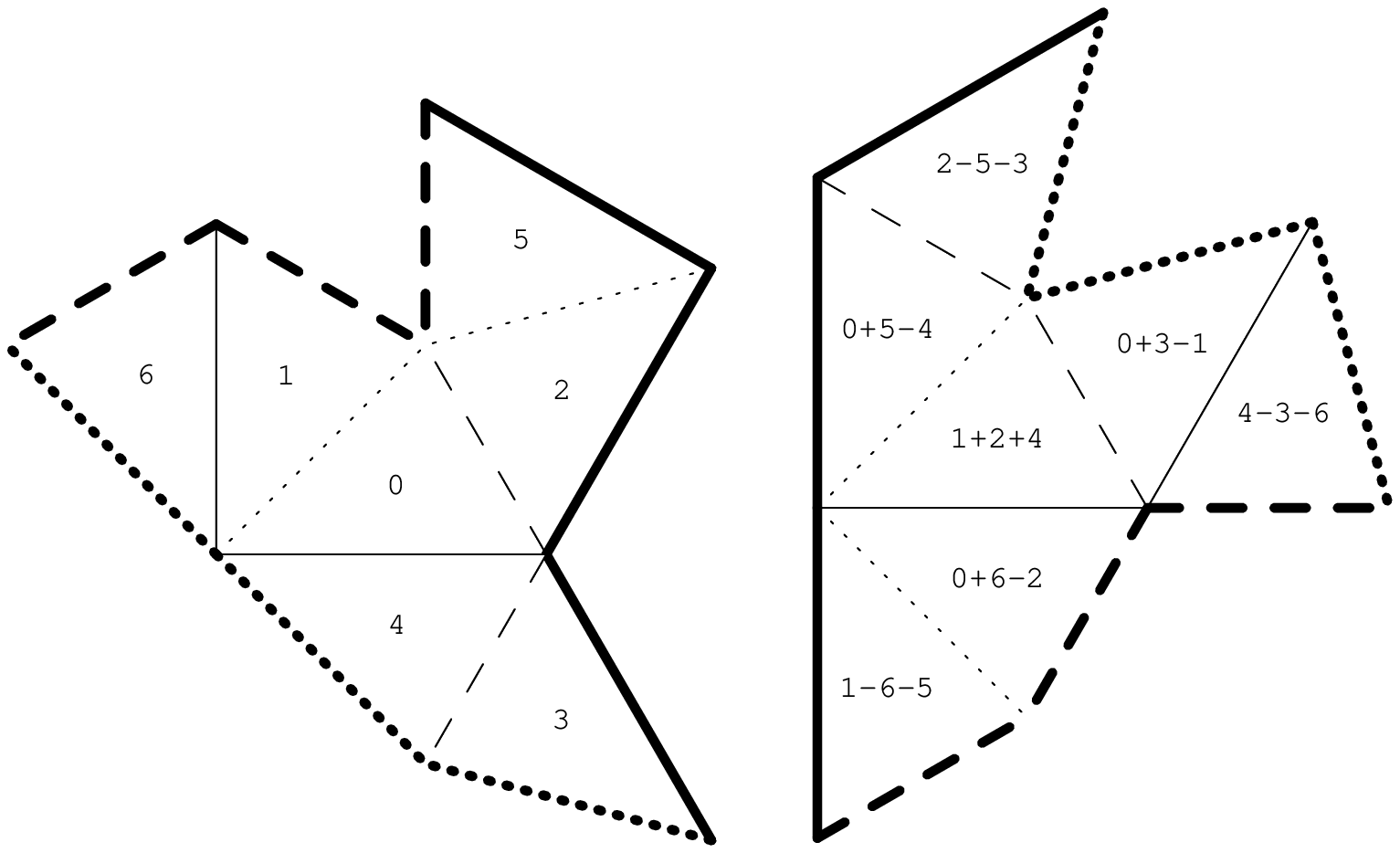}{warped}{Warped propeller.}

Now let $\lambda$  be any real number, and  $f$  any 
eigenfunction of the
Laplacian with eigenvalue  $\lambda$ for the Dirichlet problem 
corresponding
to the left-hand propeller.
Let  $f_0,f_1,\ldots,f_6$  denote the functions obtained
by restricting  $f$  to each of the  7  triangles of the left-hand 
propeller,
as indicated on the left in Figure \ref{warped}.
For brevity, we write $\0$ for $f_0$, $\1$ for $f_1$, etc.
The Dirichlet boundary condition is that  $f$  must vanish
on each boundary-segment.  Using the reflection principle, this
is equivalent to the assertion that  $f$  would
go into $-f$ if continued as a smooth eigenfunction across
any boundary-segment. (More precisely it goes into $-f \circ \sigma$
where $\sigma$ is the reflection on the boundary segment.)

On the right in Figure \ref{warped},
we show how to obtain from  $f$  another eigenfunction
of eigenvalue $\lambda$, this time for the right-hand propeller.   In 
the
central triangle, we put the function  $\1 + \2 + \4$. By this
we mean the function $f_1 \circ \tau_1 + f_2 \circ \tau_2 + f_4 \circ \tau_4$
where
for $k = 1, 2, 4$, $\tau_k$ is the isometry from the central triangle
of the right-hand propeller to the triangle labelled $k$ on the
left-hand propeller.
Now we see from the left-hand side that the functions  $\1,\2,\4$
continue smoothly across dotted lines into copies of the
functions  $\0 , \5, -\4$ respectively, so that their
sum continues into  $\0 + \5 - \4$ as shown. The reader
should check in a similar way that this continues across a solid 
line to $\4 -\5 -\0$ (its negative), and across a dashed line to
$\2 - \5 - \3$, which continues across either a solid or dotted
line to its own negative.  These
assertions, together with the similar ones obtained by symmetry (i.e. cyclic
permutation of the arms of the propellers),
are enough to show that the transplanted function is an eigenfunction
of eigenvalue  $\lambda$  that vanishes along each boundary segment 
of the right-hand propeller.

So we have defined a linear map which for each $\lambda$  takes the
$\lambda$-eigenspace for the left-hand propeller to the   
$\lambda$-eigenspace
for the right-hand one.  This is easily checked to be a
non-singular map, and so the dimension of the eigenspace
on the right-hand side is larger or equal the dimension
on the left-hand side. Since the same transplantation may also
be applied in the reversed direction the dimensions are equal.
This holds for each    $\lambda$, and so the two propellers are Dirichlet
isospectral.

In fact they are also Neumann isospectral, as can be seen by a
similar transplantation proof obtained by replacing every minus sign
in the above by a plus sign.
(Going from Neumann to Dirichlet is almost as easy:
Just color the triangles on each side alternately black and white,
and attach minus signs on the right to function elements that have moved from
black to white or vice versa.)

In the propeller example,
each of the seven function elements on the left
got transplanted into three triangles on the right, and we verified
that it all fits together seamlessly.
If we hadn't been given the transplantation rule, we could have worked it
out as follows:
We start by transplanting
the function element $\1$ into the central triangle on the right;
on the left $\1$ continues across a dotted line to $\0$, so we stick $\0$ 
in the triangle across the dotted line on the right; on the left $\0$
continues across the solid line to $\4$, and since on the right
 the solid side of the
triangle containing $\0$ is a boundary edge, we stick a
$\4$ in along with the $\0$
(don't worry about signs---we can fill them in afterwards using the
black and white coloring of the triangles);
now since on the left
$\4$ continues across a dotted line to itself
we stick a $\4$ into the center along
with the $\1$ we started with; and so on until we have three function
elements in each triangle on the right and the whole thing fits together
seamlessly.

If we had begun by putting $\0$ into
the central triangle on the right, rather than $\1$,
then we would have ended up with four function
elements in each triangle, namely,  the complement in the set
$\{\0,\1,\ldots,\6\}$ of the original three;
This gives a second transplantation mapping.
Call the original mapping $T_3$, and the complementary mapping $T_4$.
Any linear combination $a T_3 + b T_4, \; a \neq b$
will also be a transplantation mapping,
and if we take for $(a,b)$ one of the four solutions
to the equations $3 a^2 + 4 b^2 = 1$,
$a^2 + 4 a b + 2 b^2 = 0$,
our transplantation mapping
becomes {\em norm-preserving}.

Now consider the
pair of putatively homophonic domains shown in Figure \ref{homophonic} above.
In this case we find two complementary transplantation mappings
$T_5$ and $T_{16}$.
The linear combination $a T_5 + b T_{16}$ is a norm-preserving
mapping if
$5 a^2 + 16 b^2 =1$ and $a^2 + 8 a b + 12 b^2 = 0$, that is, if
$(a,b) = \pm (1/3,-1/6)$ or $(a,b) = \pm
(3/7,-1/{14})$.
In the Dirichlet case, transplantation is kind to the values of the
transplanted functions at the special interior points where six triangles
meet.
With the proper choice of sign,
the Dirichlet incarnation of $T_5$ multiplies the special value
by 2,
the Dirichlet incarnation of $T_{16}$ multiplies the special value by
$-2$,
and the four norm-preserving linear combinations
$a T_5 + b T_{16}$ specified above multiply it by
$ 2 (a-b) = \pm 1$.
Thus we can convert an orthonormal basis of Dirichlet
eigenfunctions on the left to one on the right
so that corresponding functions take on the same special value.
This shows that the two domains are homophonic,
or more specifically, {\em Dirichlet} homophonic.
There is no similar reason for these domains to be Neumann homophonic,
and, in fact, we
do not know of any pair of non-congruent Neumann homophonic domains.

\section{Gallery of examples}

\fig{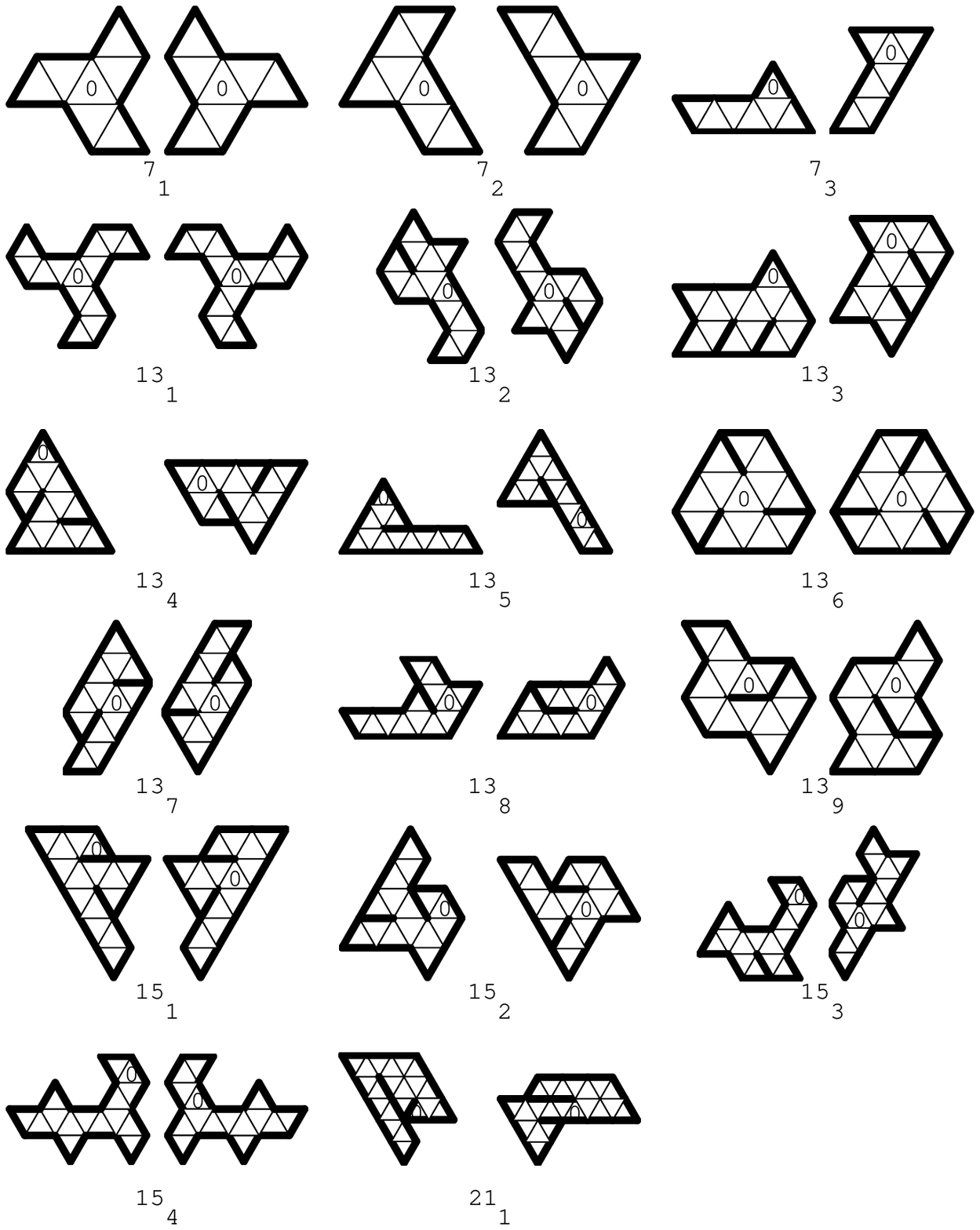}{examples}{Isospectral domains.}

Figure \ref{examples}
shows pairs of diagrams representing domains whose
isospectrality can be verified using the method of transplantation.
Each pair of diagrams represents not a single 
pair of isospectral domains, but a whole family of pairs of
isospectral domains, gotten by replacing the equilateral triangles
with general triangles so that the triangles labelled 0 are mapped
onto one another by a translation and the remaining
triangles are obtained from these by the appropriate
sequence of reflections.
We have seen two examples of this already, in
Figures \ref{warped} and \ref{homophonic}.
Further examples generated in this way are shown in Section
\ref{sec:special}.

The pair $7_1$ is the pair of propeller diagrams discussed in detail above.
The pair $7_3$ yields a simplified version of the
pair of isospectral domains
given by Gordon, Webb, and Wolpert
\cite{gordonWebbWolpert:drum},
\cite{gordonWebbWolpert:isospectral},
which was obtained by bisecting a pair of flat but non-planar isospectral
domains given earlier by Buser
\cite{buser:cayley}.
The pair $21_1$ yields the homophonic domains 
shown in Figure \ref{homophonic} above.
In this case we must be careful to choose the relevant 
angle of our generating triangle to be $2 \pi/6$ since six of these angles
meet around a vertex in each domain.
If we do not choose the angle to be $2\pi/6$, 
then instead of planar domains
we get a pair of isospectral cone-manifolds.

Note that in order for the pair $13_6$ to yield a pair
of non-overlapping non-congruent domains
we must decrease all three angles simultaneously, which we can
do by using hyperbolic triangles in place of Euclidean
triangles.  Using hyperbolic triangles, we can easily produce isospectral
pairs of {\em convex} domains in the hyperbolic plane, but we do not
know of any such pairs in the Euclidean plane.

\section{More about the examples}
The examples in the previous section were
obtained by applying a theorem of Sunada
\cite{sunada:isospectral}.
Let $G$ be a finite group.
Call two  subgroups $A$ and $B$ of $G$ {\em isospectral}
if each element of $G$ belongs to just as many conjugates of $A$
as of $B$.
(This is equivalent to requiring that $A$ and $B$
have the same number of elements in each conjugacy class of $G$.)
Sunada's theorem states that
if $G$  acts on a manifold  $M$ and $A$ and $B$ are isospectral subgroups of
$G$, then the quotient spaces of  $M$  by $A$  and  $B$  are
isospectral.

The tables in this section show for each of the examples a trio of elements
which generate the appropriate  $G$, in two distinct permutation
representations.  The isospectral subgroups  $A$  and  $B$  are the
point-stabilizers in these two permutation representations.

For the example  $7_1$, the details are as follows.
$G_0$ is the group of motions of the hyperbolic plane
generated by the reflections  $a_0,b_0,c_0$ in the sides of 
a triangle whose three angles are  $\pi/4$.
In Conway's orbifold notation
(see \cite{conway:orbifolds}),
$G_0 = *444$.
$G_0$ has a homomorphism 
$a_0 \mapsto a,b_0 \mapsto b,c_0 \mapsto c$
onto the finite group  $G=L_3(2)$ (also known as $PSL(3,2)$),
the automorphism group of the projective plane of order 2.
The generators of  $G$
act on the points and lines of this plane
(with respect to some unspecified numbering of the points and lines)
as follows:
\begin{eqnarray*}
a&=& (0\ 1)(2\ 5) \by (0\ 4)(2\ 3)\\
b&=& (0\ 2)(4\ 3) \by (0\ 1)(4\ 6)\\
c&=& (0\ 4)(1\ 6) \by (0\ 2)(1\ 5),
\end{eqnarray*}
where the actions on points and lines are separated by  $/$.

The group  $G$  has two subgroups  $A$  and  $B$  of index 7, namely
the stabilizers of a point or a line.  The preimages  $A_0$  and  $B_0$
of these two groups in  $G_0$  have fundamental regions that consist of
7 copies of the original triangle, glued together as in Figure \ref{propeller}.
Each of these is a hexagon of angles  $\pi/4,\pi/2,\pi/4,\pi/2,\pi/4,\pi/2$,
and so
each of  $A_0$  and  $B_0$  is a copy of the reflection group  $*424242$.

The preimage in $G_0$ of  the trivial subgroup of $G$ is a group
$K_0$ of index 168.  The quotient of the hyperbolic plane  by $K_0$
is a  23-fold
cross-surface (that is to say, the connected sum of  23  real projective
planes),
so that in Conway's orbifold notation
$K_0 = \cross^{23}$.
Deforming the metric on this 23-fold cross surface by replacing
its hyperbolic triangles by scalene Euclidean triangles yields a
cone-manifold $M$ whose quotients by $A$ and $B$ are non-congruent
planar isospectral domains.

Tables \ref{table:specs} and  \ref{table:perms}
display the corresponding information for our other examples.

\begin{table}
\[
\begin{array}{lllllcl}
\mbox{Pair}&\mbox{Generators}&K_0&G_0&A_0,B_0&G&\mbox{Notes}\\
\\
7_1&a,b,c&\cross^{23}&*444&*424242&L_3(2) \\
7_2&a,b',c&\cross^{16}&*443&*42423&\ditto&a' = cac\\
7_3&a',b',c&\cross^{9}&*433&*4233&\ditto&b' = aba\\
13_1&d,e,f&\cross^{704}&*444&*422422422&L_3(3)\\
13_2&d,e',f&\cross^{938}&*644&*6622342242&\ditto&e' = ded\\
13_3&d',e',f&\cross^{1172}&*664&*62234263662&\ditto&d' = fdf\\
13_4&d',e',f'&\cross^{938}&*663&*633626362&\ditto&f' = e'fe'\\
13_5&d',e'',f'&\cross^{470}&*633&*663332&\ditto&e''= d'e'd'\\
13_6&g,h,i&\cross^{1406}&*666&*632663266326&\ditto\\
13_7&g,h',i&\cross^{938}&*663&*632666233&\ditto&h' = ghg\\
13_8&g',h',i&\cross^{704}&*643&*63436222,*62633224&\ditto&g' = igi\\
13_9&g',h',i'&\cross^{938}&*644&*6262242243&\ditto&i' = g'ig'\\
15_1&j,k,l&\cross^{3362}&*663&*63362333222&L_4(2)\\
15_2&j,k,l'&\cross^{4202}&*664&*6262234342242&\ditto&l' = jlj\\
15_3&j',k,l'&\cross^{3362}&*644&*62234424242,*62422243442&\ditto&j' = kjk\\
15_4&j',k',l'&\cross^{2522}&*444&*444222442&\ditto&k' = l'kl'\\
21_1&p,q,r&\cross^{1682}&*633&*63633332,*66333323&L_3(4)\\
\end{array}
\]
\caption{Specifications.}
\label{table:specs}
\end{table}

\begin{table}
\[
\begin{array}{lll}
a&=& (0\ 1)(2\ 5) \by (0\ 4)(2\ 3)\\
b&=& (0\ 2)(4\ 3) \by (0\ 1)(4\ 6)\\
c&=& (0\ 4)(1\ 6) \by (0\ 2)(1\ 5)\\
\\
d&=&(0\ 12)(1\ 10)(3\ 5)(6\ 7) \by (0\ 4)(2\ 3)(6\ 8)(9\ 10)\\
e&=&(0\ 10)(3\ 4)(9\ 2)(5\ 8) \by (0\ 12)(6\ 9)(5\ 11)(1\ 4)\\
f&=&(0\ 4)(9\ 12)(1\ 6)(2\ 11) \by (0\ 10)(5\ 1)(2\ 7)(3\ 12)\\
\\
g&=&(0\ 2)(1\ 7)(3\ 6)(5\ 10) \by (0\ 7)(3\ 11)(6\ 8)(9\ 12)\\
h&=&(0\ 6)(3\ 8)(9\ 5)(2\ 4) \by (0\ 8)(9\ 7)(5\ 11)(1\ 10)\\
i&=&(0\ 5)(9\ 11)(1\ 2)(6\ 12) \by (0\ 11)(1\ 8)(2\ 7)(3\ 4)\\
\\
j&=&(0\ 14)(4\ 5)(9\ 10)(1\ 12)(7\ 11)(2\ 6)
\By (0\ 11)(1\ 5)(3\ 4)(6\ 10)(8\ 9)(13\ 14)\\
k&=&(4\ 6)(1\ 13)(8\ 9)(2\ 7) \By (0\ 10)(1\ 2)(6\ 9)(12\ 14)\\
l&=&(14\ 1)(3\ 4)(12\ 2)(8\ 11) \By (0\ 5)(2\ 4)(6\ 7)(11\ 14)\\
\\
p&=&(2\ 7)(3\ 11)(5\ 12)(8\ 18)(13\ 14)(15\ 17)(16\ 20) \By 
(0\ 1)(4\ 17)(7\ 12)(9\ 16)(10\ 20)(11\ 13)(15\ 19)\\
q&=&(0\ 17)(3\ 8)(4\ 12)(6\ 13)(9\ 19)(14\ 15)(16\ 18) \By 
(0\ 20)(3\ 16)(6\ 11)(8\ 15)(9\ 19)(10\ 12)(14\ 18)\\
r&=&(1\ 8)(2\ 16)(4\ 11)(5\ 19)(7\ 14)(10\ 17)(13\ 20) \By 
(1\ 8)(2\ 16)(4\ 11)(5\ 19)(7\ 14)(10\ 17)(13\ 20)

\end{array}
\]
\caption{Permutations.}
\label{table:perms}
\end{table}

Note that the permutations in Table 2 correspond to the neighboring
relations in Figure \ref{examples}. In the propeller
example, for instance, the pairs 0, 1 and 2, 5
are neighbors along a dotted line
on the left-hand side, and 0, 4 and 2, 3 are neighbors along a dotted line
on the right-hand side. Accordingly, we have the permutations 
a = (0\ 1)(2\ 5) \by (0\ 4)(2\ 3), etc.
Similar relations will hold in the other pairs of diagrams
if the triangles are properly labelled.

\section{Special cases of isospectral pairs}
\label{sec:special}

Figure \ref{special} shows some interesting special cases of isospectral pairs.

\fig{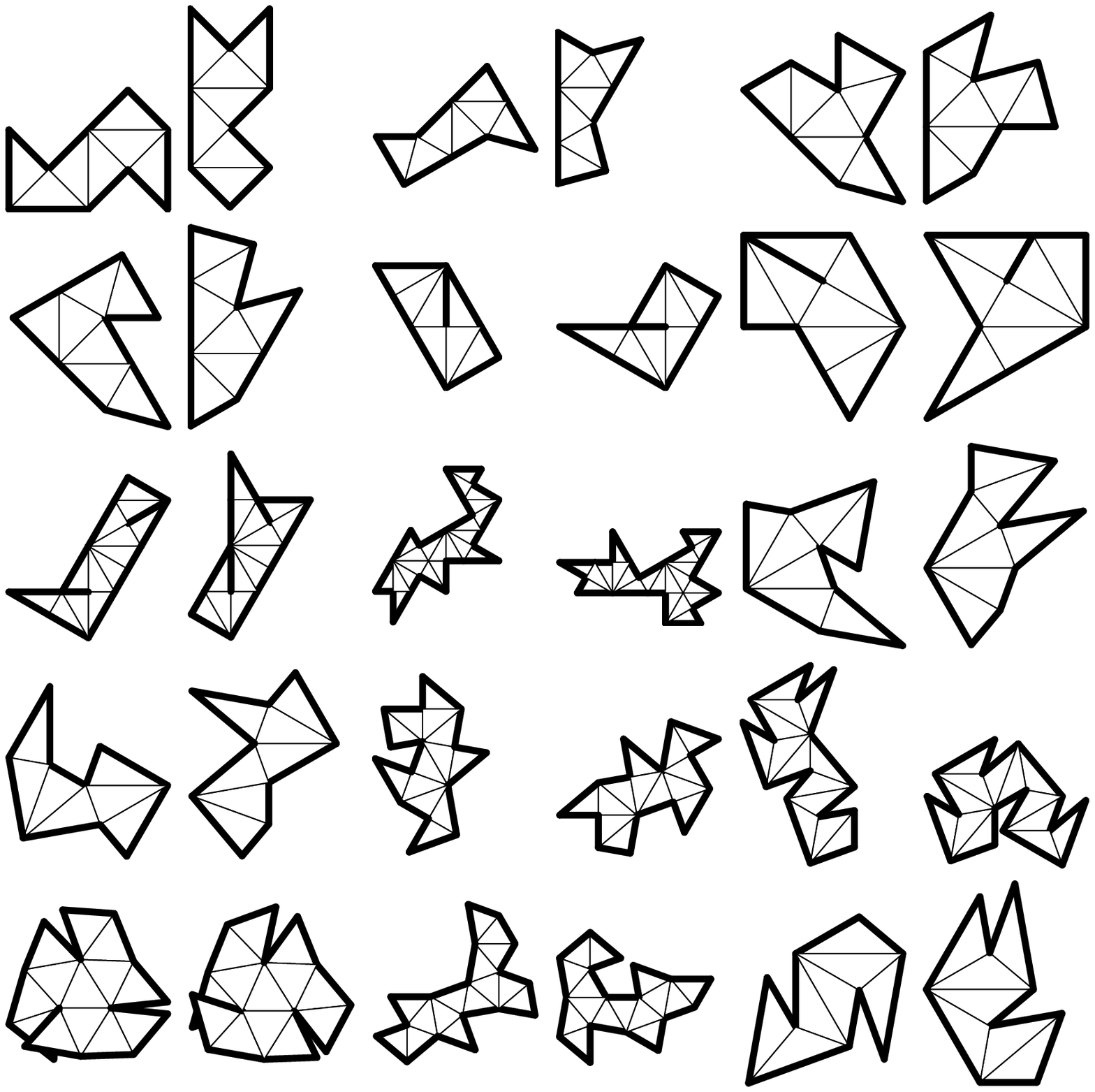}{special}{Special cases.}

\newpage
\bibliography{drum}

\begin{thebibliography}{1}

\bibitem{buser:riemann}
P.~Buser.
\newblock Isospectral {R}iemann surfaces.
\newblock {\em Ann. Inst. Fourier (Grenoble)}, 36:167--192, 1986.

\bibitem{buser:cayley}
P.~Buser.
\newblock Cayley graphs and planar isospectral domains.
\newblock In T.~Sunada, editor, {\em Geometry and Analysis on Manifolds
  (Lecture Notes in Math. 1339)}, pages 64--77. Springer, 1988.

\bibitem{conway:orbifolds}
J.~H. Conway.
\newblock The orbifold notation for surface groups.
\newblock In M.~W. Liebeck and J.~Saxl, editors, {\em Groups, Combinatorics and
  Geometry}, pages 438--447. Cambridge Univ. Press, Cambridge, 1992.

\bibitem{gordonWebbWolpert:isospectral}
C.~Gordon, D.~Webb, and S.~Wolpert.
\newblock Isospectral plane domains and surfaces via {R}iemannian orbifolds.
\newblock {\em Invent. Math.}, 110:1--22, 1992.

\bibitem{gordonWebbWolpert:drum}
C.~Gordon, D.~Webb, and S.~Wolpert.
\newblock One cannot hear the shape of a drum.
\newblock {\em Bull. Amer. Math. Soc.}, 27:134--138, 1992.

\bibitem{kac:drum}
M.~Kac.
\newblock Can one hear the shape of a drum?
\newblock {\em Amer. Math. Monthly}, 73, 1966.

\bibitem{sunada:isospectral}
T.~Sunada.
\newblock Riemannian coverings and isospectral manifolds.
\newblock {\em Ann. of Math.}, 121:169--186, 1985.

\end{thebibliography}
\bibliographystyle{plain}

\end{document}